\documentclass{amsart}
\usepackage{enumerate, url}
\usepackage{amssymb}

\usepackage{color}
\usepackage{graphicx}
\usepackage{pgfplots}
\usepackage{xcolor}

\usepackage{amsmath}
\usepackage{enumerate}
\usepackage{bbm}
\usepackage{verbatim}
\usepackage{amsthm}
\usepackage{bbm}

\chardef\bslash=`\\ 

\newtheorem{thm}{Theorem}[section]

\newtheorem*{openproblem}{Open Problem}

\theoremstyle{definition}

\newtheorem{example}{Example}

\theoremstyle{remark}

\newcommand{\NN}{\mathbb{N}} 
\newcommand{\Span}{\mathrm{Span}} 
\newcommand*\abs[1]{\lvert#1\rvert} 
\newcommand{\Sym}{\mathrm{Sym}}
\newcommand{\SoS}{\mathrm{SoS}}
\newcommand{\ii}{\mathrm{i}}

\newcommand{\eval}[2][\right]{\relax
	\ifx#1\right\relax \left.\fi#2#1\rvert}

\let\abs=\envert

\begin{document}
	\title[A Negative Answer]{The non-Archimedean  Nirgendsnegativsemidefinitheitsstellensatz  is not true
		\thanks{This research is supported by the National Key Research
			Project of China 2018YFA0306702  and the National Natural Science
			Foundation of China 12071467}
	}
	\author[H. Liang, S. Yan, J. Yang, L. Zhi]{Hao Liang, Sizhuo Yan, Jianting Yang, Lihong Zhi}
	\address{ KLMM, Academy of Mathematics and Systems Science, Chinese
		Academy of Sciences,
		University of Chinese Academy of Sciences, Beijing 100190, China
	}
	
	\email{\{lianghao2020,yansizhuo,yangjianting\}@amss.ac.cn, lzhi@mmrc.iss.ac.cn}

	
	%
	%
	%
	%

	\begin{abstract}
		Klep and Schweighofer asked whether the Nirgendsnegativsemide-finitheitsstellensatz  holds for a symmetric noncommutative polynomial whose evaluations at bounded self-adjoint operators on any nontrivial  Hilbert space are not negative semidefinite.
		We provide an example to show the open problem has a negative answer.
	\end{abstract}
	
	\maketitle
	\tableofcontents

	\section{Introduction}

	This paper considers  polynomials generated by  noncommutative variables $\underline{X}:=\{X_1, X_2, \ldots, X_m\}$ with coefficients from  $k \in\{\mathbb{C},\mathbb{R}\}$,  where $ \mathbb{R}$, $ \mathbb{C}$ are real numbers and complex numbers respectively. Let  $\mathbb{N}:=\{1,2,\ldots\}$ be the set of natural numbers. Let
	$$p=\sum_{\omega \in \mathcal{W}_m} p_{\omega} \omega$$
	be a polynomial in   $k\left\langle \underline{X}\right\rangle $  with
	finitely many nonzero coefficients $p_\omega$, and let $\mathcal{W}_{m}$ be a set of words generated by  $\{X_1,X_2, \ldots,X_m\}$. The length of the longest word appearing in a polynomial $p$ is defined as  the degree of $p$. We define the transpose of a polynomial $p$ as
	$$
	p^{*}=\sum_{\omega \in \mathcal{W}_m} p_{\omega}^{*} \omega ^{*},
	$$
	where $\omega^*=X_{i_k}\cdots X_{i_2} X_{i_1}$
	for the word $\omega=X_{i_1}X_{i_2}\cdots X_{i_k}$.
	If $p=p^{*}$, we say $p$ is symmetric. The set of symmetric polynomials is denoted by  $\Sym\,k\langle \underline{X}\rangle$.


	Let  $\mathcal{H}$ denote a separable  $k-$Hilbert space,  and let $B(\mathcal{H})$ denote the set of  bounded operators on  $\mathcal{H}$.
	We  evaluate a polynomial $p$ at $\underline{A}=(A_1,A_2,\ldots,A_m)$, where each $A_i$ is a self-adjoint operator in $B(\mathcal{H})$ for  $1\leq i \leq m$.
	The evaluation of a tuple $\underline{A}=(A_1,A_2,\ldots,A_m)$ on the empty word is $\rm{Id}_\mathcal{H}$, which is the identity
	operator in $\mathcal{H}$.
	
	Let $S$ be  a subset of symmetric polynomials $\Sym\,k\langle \underline{X}\rangle$. The  noncommutative semialgebraic set  $K_S$  consists  of  tuples $\underline{A}=(A_1, \ldots, A_m)$
	of  bounded self-adjoint operators on a nontrivial $k$-Hilbert space $\mathcal{H}$   such that $s(\underline{A}) $ is positive semidefinite for all $ s \in S$, i.e.,
	\[K_S=\left\{\underline{A} \in B(\mathcal{H})^m~|~A_i^*=A_i,~0\leq i\leq m, ~s(\underline{A}) \succeq 0 ,~\forall s \in S\right\}.\]
	
	Let $M_S$  be  the quadratic module defined by  the set of  elements of the form
	\[M_S=\left \{\sum_{i=1}^{n} g_i^* s_i g_i~|~g_i\in k\left\langle \underline{X}\right\rangle,~s_i\in S \cup \{1\},~n\in \mathbb{N}\right\}.\]
	When $S$ is an empty set, the quadratic module $M_\emptyset$ is denoted by  $\SoS$.
	
	If $M_S$ is an Archimedean quadratic module, i.e., there exists an
	$N\in \mathbb{N}$, such that $N-X_1^2-X_2^2-\cdots-X^2_m \in M_S$, then it has been proved in  Theorem 1.4  \cite{IgorMarkus2007nicht} that
	$p(\underline{A})$ is not negative semidefinite for any $\underline{A} \in K_S$ 
	if and only if there exist a positive integer $r \in \mathbb{N} $ and polynomials $g_1,\dots g_r\in k\langle\underline{X}\rangle$ such that $$\sum_{i=1}^r{g_i^* pg_i}\in 1+M_S.$$
	
	Different proofs have  been given in
	Theorem 5 in  \cite{cimprivc2008maximal},  Theorem 4.2 in \cite{cimprivc2011noncommutative}, and  Proposition 17 in \cite{schmudgen2009noncommutative}.
	
	
	\section{A negative answer to  Klep and Schweighofer's open problem}
	Klep and Schweighofer  posted an interesting open problem in \cite{IgorMarkus2007nicht}.
	
	\begin{openproblem}\label{openproblem}
		\cite{IgorMarkus2007nicht} Given a symmetric polynomial  $p\in  k\langle\underline{X}\rangle$, are the following two conditions  equivalent?
		\begin{enumerate}
			\item[$(a)$] $p(A_1, \ldots , A_m)$ is not negative semidefinite for any  nontrivial $k$-Hilbert space $\mathcal{H}$ ($\mathcal{H}\neq0$) and bounded self-adjoint operators $A_1,\dots A_m$ on $\mathcal{H}$;
			\item[$(b)$] There exist $r \in \mathbb{N} $ and  polynomials $g_1,\dots g_r\in k\langle\underline{X}\rangle$  such that
			\begin{equation}\label{con2}
				\sum_{i=1}^r{g_i^* p g_i}\in1+{\rm SoS}.
			\end{equation}
		\end{enumerate}
	\end{openproblem}
	If $f$ is a commutative polynomial in $ \mathbb{R}[X_1, \ldots, X_m]$, then we have
	\[f(A_1, \ldots , A_m) \not\preceq 0 \iff f(A_1, \ldots , A_m) > 0, ~{\rm where}~ A_i \in \mathbb{R}, 1\leq i \leq m.\]
	According to  the Positivstellensatz  in \cite{ krivine1964anneaux, stengle1974nullstellensatz}, we have $f>0$ on $\mathbb{R}^m$ if and only if there exist $r \in \mathbb{N} $ and polynomials $g_1,\dots g_r\in \mathbb{R}[X_1, \ldots, X_m]$ such that
	\[ \left( \sum_{i=1}^r g_i^2 \right) f \in 1+{\rm SoS}.  \]
	Therefore, for commutative cases ($\dim \mathcal{H} =1$), Klep and Schweighofer's question has a positive answer.  However, their open problem has a negative answer   for symmetric  noncommutative polynomials whose evaluations at bounded operators on  any nontrivial Hilbert space  are not negative semidefinite.

	
	\begin{example}\label{example1} Let the polynomial $f(X_1,X_2)$  be  given as
		\begin{equation}\label{ex1}	
			f(X_1,X_2)=X_1 X_2^2 X_1-X_2 X_1^2 X_2+1.
		\end{equation}
	\end{example}
	
	\begin{thm}\label{thm2.1}

		The polynomial $f$ given in (\ref{ex1}) satisfies condition $(a)$ but does not satisfy  condition $(b)$ in the open problem. Hence, the non-Archimedean  Nirgendsnegativsemidefinitheitsstellensatz  is not true.
		
	\end{thm}
	
	\begin{proof}
		Firstly, we show that $f$  satisfies  condition $(a)$,  $$ f(A,  B) \not\preceq 0, ~\forall A=A^*, B=B^*,  A, B \in  B(\mathcal{H}).$$
		
		For any Hilbert space $\mathcal{H}$ and bounded self-adjoint operator $A, B \in B(\mathcal{H})$, we have $AB=(BA)^*$. Therefore,  the bounded operator $AB^2A$ and $BA^2B$ are both positive semidefinite. Moreover,   the set of  bounded operators $B(\mathcal{H})$ is an example of $C^*$-algebra \cite{Arveson1976AnIT}. Let
		the symbol $\|\cdot\|$ denote the operator norm.  We have the following equations
		\[a=\|AB^2A\|=\|(BA)^*BA\|=\|BA\|^2=\|A B\|^2= \|BA^2B\|.\]
		
		If $a=0$, then we have $AB^2A=BA^2B=0$, and $f(A,B)={\rm Id_{\mathcal{H}}}$ is not negative semidefinite.
		
		If $a>0$, let  $\langle \cdot,\cdot \rangle$ denote inner product on $\mathcal{H}$.  There exists  a sequence of unit vectors $\{v_i\}\in \mathcal{H}$ such that
		\[a=\lim\limits_{i\rightarrow \infty}\langle AB^2Av_i,   v_i \rangle.\]
		Hence, there exists a number $j\in \mathbb{N}$ such that
		\[\langle AB^2Av_j, v_j \rangle>a-\frac{1}{2},\quad -\langle BA^2Bv_j, v_j \rangle \geq -a.\]
		Adding two sides of the  inequalities,  we derive that  $f(A,B)$ is not negative semidefinite as
		\[ \langle f(A,B)v_j, v_j \rangle >\frac{1}{2}.\]

		Secondly, we show that  $f$ does not  satisfy  condition $(b)$.
		In \cite[Lemma 2.1]{Helton2002positive}, Helton proved that for any given  symmetric noncommutative polynomial $p$, there exists a symmetric matrix $\mathcal{M}_p$, not dependent on $\underline{X}$, and a vector $V(\underline{X})$ of monomials in $\underline{X}$ such that
		\[p(\underline{X})=V(\underline{X})^T \mathcal{M}_p V(\underline{X}).\]
		Furthermore, the vector $V(\underline{X})$ can always be chosen as $V^d(\underline{X})$ which  contains the monomials whose degree is at most $d$,  where $d=\lceil \deg(p)/2 \rceil$.  As shown in \cite[Theorem 2.1]{mccullough2005noncommutative},
		a symmetric polynomial $p\in \SoS$ if and only if there exists a positive semidefinite matrix $\mathcal{M}_p$ such that
		\[p(\underline{X})=V(\underline{X})^T \mathcal{M}_p V(\underline{X}).\]
		
		For a monomial $\omega$ in $p$, we can find its corresponding entries whose index $(v,u)$ satisfies $\omega=v^*u$ in $\mathcal{M}_p$, and the coefficient of $\omega$ in $p$ equals to $\sum_{\omega=v^*u}\mathcal{M}_p(v,u)$.    The matrix $\mathcal{M}_p$ is not unique, but for a monomial $\omega$ in $p$ whose degree is $2d$, it has only one corresponding entry, i.e., there is a unique choice such that $\omega=v^*u$ where degrees of $u$ and $v$ are equal to $d$.
		
		For the polynomial $f$ defined in   (\ref{ex1}), we have  
		\[f(X_1,X_2)=V^2(\underline{X})^* \mathcal{M}_f V^2(\underline{X})\]
		where a matrix $\mathcal{M}_f$ can be  written  as
		\[\mathcal{M}_f= \bordermatrix{  &1 &X_1 &X_2 &X_1^2 &X_1 X_2 &X_2 X_1 &X_2^2\cr
			1 &1 & & & & & &\cr
			X_1 & & 0& & & & &\cr
			X_2& & & 0& & && \cr
			X_1^2& & & &0 & & & \cr
			X_1 X_2 & & & & & -1 & &\cr
			X_2X_1 & & & & & & 1& \cr
			X_2^2 & & & & & & &0}
		\]
		and $V^2(\underline{X})$ is the monomial vector satisfying
		\[
		V^2(\underline{X})^*=\begin{pmatrix}
			1&X_1&X_2&X_1^2&X_2X_1&X_1 X_2&X_2^2
		\end{pmatrix}.\]

		
		For any $r\in \mathbb{N}$ and $g_1,g_2,\ldots,g_r \in k\langle \underline{X} \rangle$, we have the symmetric polynomial
		\[F=\sum_{i=1}^{r}g_i^* f g_i=\sum_{i=1}^{r}g_i^* X_1 X_2^2 X_1 g_i -\sum_{i=1}^{r}g_i^* X_2X_1^2 X_2 g_i +\sum_{i=1}^{r}g_i^* g_i.\]
		Define $D=\max\{d_1,\ldots,d_r\}$, where $d_i$ denotes the degree of $g_i$ for  $1\leq i \leq r$.
		There always exists  a monomial $u$   with the maximal degree $D$    in some $g_i$  for   $ 1 \leq i \leq r$. 
		
		For any matrix $\mathcal{M}_0$ which satisfies that
		\[F=V^{D+2}(\underline{X}) \mathcal{M}_0 V^{D+2}(\underline{X}),\]
		the diagonal entry in $\mathcal{M}_0$ indexed by $(X_1X_2u,X_1X_2u)$ must be the coefficient of the monomial $u^*X_2X_1^2X_2u$ which is negative in $F$  since the degree of $u^*X_2X_1^2X_2u$ is $2D+4$. Therefore,  $\mathcal{M}_0$ can not be positive semidefinite, i.e., we have
		$$F=\sum_{i=1}^{r}g_i^* f g_i\notin \SoS.$$
		Since  $1+{\rm SoS}\subseteq \SoS$, we have shown that  there are no $r\in \mathbb{N}$ and $g_1,g_2,\ldots,g_r \in k\langle \underline{X} \rangle$, such that  condition  $(b)$ holds.
	\end{proof}
	
	\section{Further discussion}
	
	We give below  a symmetric polynomial $g(X_1,X_2)$  with complex coefficients
	which satisfies condition $(a)$ but does not satisfy  condition $(b)$ in the open problem.   It is worth noting that,   if we allow evaluations at unbounded operators, then the polynomial in  Example \ref{ex3}  does not satisfy condition $(a)$.

	\begin{example}\label{ex3}
		Let the polynomial $g(X_1,X_2)$ be given as 
		
		\begin{equation}\label{example3}
			g(X_1,X_2)=2\ii  X_2X_1-2\ii  X_1X_2+1.
		\end{equation}
	\end{example}
	

	\begin{thm}\label{thm3.1}
		The symmetric polynomial $g$ given in (\ref{example3}) satisfies condition $(a)$ but does not satisfy  condition $(b)$ in the open problem. 
	\end{thm}


	The polynomial $g(X_1,X_2)$ given in (\ref{example3}) actually is obtained by replacing  the non-symmetric variable $X$ 
	in the following polynomial $h(X,X^*)$:
	\begin{equation}\label{ex2}
		h(X,X^*)=XX^*-X^*X+1
	\end{equation}
	by two symmetric variables
	\begin{equation}
		X_1=\frac{X+X^*}{2}, X_2=\frac{X-X^*}{2\ii}. 
	\end{equation}
	We have 
	\begin{equation}\label{relation}
		g(X_1,X_2)=g\left(\frac{X+X^*}{2},\frac{X-X^*}{2\ii}\right)=h(X,X^*).
	\end{equation}
	We use the  polynomial $h(X,X^*)$ and the relation (\ref{relation}) between $h(X,X^*)$ and $g(X_1,X_2)$  to simplify the   proof of Theorem \ref{thm3.1}.  
	
	\begin{proof}
		
		Similar to the proof in Theorem \ref{thm2.1}, we can show that 
		$h(A,A^*)$ is not negative semidefinite for any  nontrivial $k$-Hilbert space $\mathcal{H}$ ($\mathcal{H}\neq0$) and bounded operators $A, A^*\in B(\mathcal{H})$, i.e., 
		\begin{equation}\label{hconda}
			h(A,  A^*) \not\preceq 0, ~\forall  A, A^*  \in  B(\mathcal{H}).
		\end{equation}

		According to (\ref{relation}) and (\ref{hconda}), we know that  $g(X_1,X_2)$ satisfies  condition $(a)$, i.e.,
		\[ g(A,  B) \not\preceq 0, ~\forall A=A^*, B=B^*,  A, B \in  B(\mathcal{H}).\]

		Now let us show 
		that the polynomial $g(X_1,X_2)$  does not satisfy  condition $(b)$. 
		The proof is similar to the one given in Theorem \ref{thm2.1}.
		For any $r\in \mathbb{N}$, and $g_1,g_2,\ldots,g_r \in k\langle \underline{X} \rangle$, we have the symmetric polynomial
		\[G=\sum_{j=1}^{r}g_j^* g  g_j=\sum_{j=1}^{r}2\ii g_j^* X_2 X_1 g_j -\sum_{j=1}^{r}2\ii g_j^* X_1X_2 g_j +\sum_{j=1}^{r}g_j^* g_j.\]
		Define $D=\max\{d_1,\ldots,d_r\}$, where $d_j$ denotes the degree of $g_j$ for all $1\leq j \leq r$. 
		There always exists  a monomial $u$   with the maximal degree $D$    in some $g_j$  for   $ 1 \leq j \leq r$,
		and  the degrees of $u^*X_2X_1u$ and $u^*X_1X_2u$ are $2D+2$.
		
		For any matrix $\mathcal{M}_0$ which satisfies that
		\[G=V^{D+1}(\underline{X}) \mathcal{M}_0 V^{D+1}(\underline{X}),\]
		the principal minor in $\mathcal{M}_0$ indexed by $X_1u$ and $X_2u$ must have the following form:
		\[\bordermatrix{  &X_1u &X_2u \cr
			X_1u & 0& -2\ii a\cr
			X_2u&2\ii a & 0\cr}, \quad a \neq 0\in \mathbb{R}.
		\]
		The determinant of the minor   is $-4a^2<0$. 
		Therefore,  $\mathcal{M}_0$ can not be positive semidefinite, i.e., we have
		$$G=\sum_{j=1}^{r}g_j^* g g_j\notin 1+ \SoS.$$
	\end{proof}

	It is interesting to notice that if we allow evaluations at unbounded operators \cite{schmudgen2020invitation}, the polynomial $g$ given in (\ref{example3}) does not satisfy  condition $(a)$.

	Let us  consider below  unbounded operators $(A, D(A))$ on $\mathcal{H}$ with a dense domain $D(A)$ such that
	\begin{equation}\label{unbound}
		A(D(A))\subseteq D(A), ~D(A)\subseteq D(A^*), ~A^*(D(A))\subseteq D(A).
	\end{equation}
	We generalize  the positivity  by the following condition \cite{schmudgen2009noncommutative}:
	\[p(A,A^*) \succeq 0 ~{\rm if ~ and ~ only ~if}~  \langle p(A,A^*) \phi, \phi \rangle \geq 0 ~{\rm for ~all}~ \phi \in D(A). \]
	
	Take an orthonormal basis $(e_k:k\in\NN)$ of $\mathcal{H}$. Define
	\[A e_k=ke_{k+1} ~{\rm for}~k\in\NN,\]
	and let
	\[D(A)=\Span(e_k:k\in\NN)\]
	be the minimal linear subspace of $\mathcal{H}$ containing $(e_k:k\in\NN)$. Thus, the matrix of $A$ under the basis $(e_k:k\in\NN)$ is
	\[\begin{pmatrix}
		0& & & &\\
		1&0& & &\\
		&2&0&&\\
		& &\ddots&\ddots
	\end{pmatrix}.\]
	For any vector $v=\sum_{k=1}^{\infty}{v_k e_k}\in D(A)$, we have
	$$\langle  (h(A,A^*) v,v \rangle = \langle  (AA^*-A^*A+{\rm Id_\mathcal{H}})v,v \rangle =-2\sum_{k=1}^{\infty}k\abs{v_{k+1}}^2\leq 0.$$
	
	Let 
	\[B=\frac{A+A^*}{2}, ~C=\frac{A-A^*}{2\ii}. \]
	According to (\ref{relation}),
	for any vector $v=\sum_{k=1}^{\infty}{v_k e_k}\in D(A)$, we have
	\[\langle  g(B,C) v,v \rangle=  \langle  (2\ii C B-2\ii B C+{\rm Id_\mathcal{H}})v,v \rangle =-2\sum_{k=1}^{\infty}k\abs{v_{k+1}}^2\leq 0.\]
	Hence, we have proved that $g(B, C)$ is negative semidefinite for the unbounded operators $(B,D(B))$ and $(C,D(C))$. Hence the polynomial  $g$ does not satisfy  condition $(a)$ in the open problem.

	We know that  condition $(b)$ is stronger than condition $(a)$ in the open problem. 
	Inspired by  Schm{\"u}dgen's survey paper \cite{schmudgen2009noncommutative} and his books on unbounded operators  and representations \cite{schmudgen1990unbounded,schmudgen2020invitation}, we replace the bounded operators with operators satisfying (\ref{unbound}) in condition $(a)$. It is clear that condition $(a')$ is stronger than condition $(a)$. We  ask the following question:

	\begin{openproblem}\label{openproblem2}
		Given a symmetric polynomial  $p\in  k\langle\underline{X}\rangle$, are the following two conditions  equivalent?
		\begin{enumerate}
			\item[$(a')$] $p(A_1, \ldots , A_m)$ is not negative semidefinite for any  nontrivial $k$-Hilbert space $\mathcal{H}$ ($\mathcal{H}\neq0$) and self-adjoint operators $A_1,\dots A_m$ on $\mathcal{H}$ satisfying (\ref{unbound});
			\item[$(b)$] There exist $r \in \mathbb{N} $ and    polynomials $g_1,\dots g_r\in k\langle\underline{X}\rangle$ such that
			\begin{equation*}\label{2con2}
				\sum_{i=1}^r{g_i^* pg_i}\in1+{\rm SoS}.
			\end{equation*}
		\end{enumerate}
	\end{openproblem}

	\noindent{\bfseries Acknowledgments:}
	We are grateful to the anonymous referee for the constructive
	comments. 
	We would  like to thank  Zhihong Yang and Tianshi Yu for  helpful
	discussions. This research is supported by the National Key Research
	Project of China 2018YFA0306702  and the National Natural Science
	Foundation of China 12071467.
	
	\bibliography{ce}
	\bibliographystyle{ijmart}
	
\end{document}